\documentclass[11pt]{article}
\usepackage{mathrsfs}
\usepackage{amscd}
\usepackage{amsmath,amsfonts,amssymb,amscd}
\usepackage{indentfirst,graphics,epsfig,psfrag}
\input{epsf}
\usepackage{ifpdf}
\usepackage{enumerate}
\usepackage{appendix}
\usepackage{enumerate}
\usepackage{lineno}

\makeatletter \@addtoreset{figure}{section} \makeatother
\makeatletter
\long\def\@makecaption#1#2{%
   \vskip 10\p@
   \setbox\@tempboxa\hbox{{#1}\ \ #2}%
   \ifdim \wd\@tempboxa >\hsize

       {#1}\ \ #2\par
   \else
       \hbox to\hsize{\hfil\box\@tempboxa\hfil}%
   \fi}
\makeatother

\newtheorem{thm}{Theorem}[section]
\newtheorem{cor}[thm]{Corollary}
\newtheorem{lem}[thm]{Lemma}

\newtheorem{obs}[thm]{Observation}
\newtheorem{pro}[thm]{Proposition}

\newcommand{\qed}{{\hfill\rule{3pt}{7pt}}}

\def\qed{\hfill \rule{4pt}{7pt}}

\begin{document}
\title{\textbf{Path connectivity
of line graphs}\footnote{Supported by the National Science
Foundation of China (Nos. 11551001, 11161037, 11461054) and the
Science Found of Qinghai Province (No. 2014-ZJ-907).}}
\author{
\small Yaping Mao$^{1,2}$\footnote{E-mail: maoyaping@ymail.com}\\[0.2cm]
\small $^{1}$Department of Mathematics, Qinghai Normal University,\\[0.1mm]
\small $^{2}$Key Laboratory of IOT of Qinghai Province,\\[0.1mm]
\small Xining, Qinghai 810008, China}
\date{}
\maketitle
\begin{abstract}
Dirac showed that in a $(k-1)$-connected graph there is a path
through each $k$ vertices. The path $k$-connectivity $\pi_k(G)$ of a
graph $G$, which is a generalization of Dirac's notion, was
introduced by Hager in 1986. In this
paper, we study path connectivity of line graphs.\\[2mm]
\end{abstract}

\section{Introduction}

An interpersonal network is represented as a graph, where a node is
a processor and an edge is a communication link. Broadcasting is the
process of sending a message from the source node to all other nodes
in a network. It can be accomplished by message dissemination in
such a way that each node repeatedly receives and forwards messages.
Some of the nodes and/or links may be faulty. However, multiple
copies of messages can be disseminated through disjoint paths. We
say that the broadcasting succeeds if all the healthy nodes in the
network finally obtain the correct message from the source node
within a certain limit of time. A lot of attention has been devoted
to fault-tolerant broadcasting in networks \cite{Fragopoulou,
Hedetniemi, Jalote, Ramanathan}. In order to measure the degree of
fault-tolerance, the above disjoint path structure connecting two
nodes is generalized into some tree structures connecting more than
two nodes, see \cite{Ku, LLSun, LM3}. To show the properties of
these generalizations clearly, we hope to state from the
connectivity in Graph Theory. We divide our introduction into the
following four subsections to state the motivations and our results
of this paper.

\subsection{Connectivity and $k$-connectivity}

All graphs considered in this paper are undirected, finite and
simple. We refer to \cite{bondy} for graph theoretical notation and
terminology not described here. For a graph $G$, let $V(G)$, $E(G)$
and $\delta(G)$ denote the set of vertices, the set of edges and the
minimum degree of $G$, respectively. For $S\subseteq V(G)$, we
denote by $G-S$ the subgraph obtained by deleting from $G$ the
vertices of $S$ together with the edges incident with them.
Connectivity is one of the most basic concepts in graph theory, both
in combinatorial sense and in algorithmic sense. As we know, the
classical connectivity has two equivalent definitions. The
\emph{connectivity} of $G$, written $\kappa(G)$, is the minimum size
of a vertex set $S\subseteq V(G)$ such that $G-S$ is disconnected or
has only one vertex. We call this definition the `cut' version
definition of connectivity. A well-known theorem of Whitney
\cite{Whitney} provides an equivalent definition of connectivity,
which can be called the `path' version definition of connectivity.
For any two distinct vertices $x$ and $y$ in $G$, the \emph{local
connectivity} $\kappa_{G}(x,y)$ is the maximum number of internally
disjoint paths connecting $x$ and $y$. Then
$\kappa(G)=\min\{\kappa_{G}(x,y)\,|\,x,y\in V(G),x\neq y\}$ is
defined to be the \emph{connectivity} of $G$. Similarly, the
classical edge-connectivity also has two equivalent definitions. The
\emph{edge-connectivity} of $G$, written $\lambda(G)$, is the
minimum size of an edge set $M\subseteq E(G)$ such that $G-M$ is
disconnected or has only one vertex. We call this definition the
`cut' version definition of edge-connectivity. For any two distinct
vertices $x$ and $y$ in $G$, the \emph{local edge-connectivity}
$\lambda_{G}(x,y)$ is the maximum number of edge-disjoint paths
connecting $x$ and $y$. Then
$\lambda(G)=\min\{\lambda_{G}(x,y)\,|\,x,y\in V(G),x\neq y\}$ is
defined to be the \emph{edge-connectivity} of $G$, which can be
called the `path' version definition; see \cite{FordF}. For
connectivity and edge-connectivity, Oellermann published a survey
paper; see \cite{Oellermann2}.

Although there are many elegant and powerful results on connectivity
in graph theory, the classical connectivity and edge-connectivity
cannot be satisfied considerably in practical uses. So people want
some generalizations of both connectivity and edge-connectivity. For
the `cut' version definition of connectivity, we are looking for a
minimum vertex-cut with no consideration about the number of
components of $G-S$. Two graphs with the same connectivity may have
different degrees of vulnerability in the sense that the deletion of
a vertex cut-set of minimum cardinality from one graph may produce a
graph with considerably more components than in the case of the
other graph. For example, the star $K_{1,n}$ and the path $P_{n+1}\
(n\geq 3)$ are both trees of order $n+1$ and therefore connectivity
$1$, but the deletion of a cut-vertex from $K_{1,n}$ produces a
graph with $n$ components while the deletion of a cut-vertex from
$P_{n+1}$ produces only two components. Chartrand et al.
\cite{Chartrand1} generalized the `cut' version definition of
connectivity. For an integer $k \ (k\geq 2)$ and a graph $G$ of
order $n \ (n\geq k)$, the \emph{$k$-connectivity} $\kappa'_k(G)$ is
the smallest number of vertices whose removal from $G$ produces a
graph with at least $k$ components or a graph with fewer than $k$
vertices. Thus, for $k=2$, $\kappa'_2(G)=\kappa(G)$. For more
details about $k$-connectivity, we refer to \cite{Chartrand1, Day,
Oellermann2, Oellermann3}. The $k$-edge-connectivity, which is a
generalization of the `cut' version definition of classical
edge-connectivity was initially introduced by Boesch and Chen
\cite{Boesch} and subsequently studied by Goldsmith in
\cite{Goldsmith1, Goldsmith2} and Goldsmith et al.
\cite{Goldsmith3}. For more details, we refer to \cite{Beineke,
Oellermann1}.

\subsection{Generalized connectivity and generalized edge-connectivity}

The generalized connectivity of a graph $G$, introduced by Hager
\cite{Hager}, is a natural and nice generalization of the `path'
version definition of connectivity. For a graph $G=(V,E)$ and a set
$S\subseteq V$ of at least two vertices, \emph{an $S$-Steiner tree}
or \emph{a Steiner tree connecting $S$} (or simply, \emph{an
$S$-tree}) is a subgraph $T=(V',E')$ of $G$ that is a tree with
$S\subseteq V'$. Two Steiner trees $T$ and $T'$ connecting $S$ are
said to be \emph{internally disjoint} if $E(T)\cap
E(T')=\varnothing$ and $V(T)\cap V(T')=S$. For $S\subseteq V(G)$ and
$|S|\geq 2$, the \emph{generalized local connectivity} $\kappa(S)$
is the maximum number of internally disjoint Steiner trees
connecting $S$ in $G$. Note that when $|S|=2$ a minimal Steiner tree
connecting $S$ is just a path connecting the two vertices of $S$.
For an integer $k$ with $2\leq k\leq n$, \emph{generalized
$k$-connectivity} (or \emph{$k$-tree-connectivity}) is defined as
$\kappa_k(G)=\min\{\kappa(S)\,|\,S\subseteq V(G),|S|=k\}$. Clearly,
when $|S|=2$, $\kappa_2(G)$ is nothing new but the connectivity
$\kappa(G)$ of $G$, that is, $\kappa_2(G)=\kappa(G)$, which is the
reason why one addresses $\kappa_k(G)$ as the generalized
connectivity of $G$. By convention, for a connected graph $G$ with
less than $k$ vertices, we set $\kappa_k(G)=1$. Set $\kappa_k(G)=0$
when $G$ is disconnected. This concept appears to have been
introduced by Hager in \cite{Hager}. It is also studied in
\cite{Chartrand2} for example, where the exact value of the
generalized $k$-connectivity of complete graphs are obtained. Note
that the generalized $k$-connectivity and the $k$-connectivity of a
graph are indeed different. Take for example, the graph $H_1$
obtained from a triangle with vertex set $\{v_1,v_2,v_3\}$ by adding
three new vertices $u_1,u_2,u_3$ and joining $v_i$ to $u_i$ by an
edge for $1 \leq i\leq 3$. Then $\kappa_3(H_1)=1$ but
$\kappa'_3(H_1)=2$. There are many results on the generalized
connectivity or tree-connectivity, we refer to \cite{Chartrand2,
LLM, LLSun, LM1, LM2, LM3, LM4, LMS, Okamoto}.

The following Table 1 shows how the generalization proceeds. {\small
\begin{center}
\begin{tabular}{|c|c|c|}
\hline  & Classical connectivity& Generalized
connectivity\\[0.1cm]
\cline{1-3}
Vertex subset & $S=\{x,y\}\subseteq V(G) \ (|S|=2)$ & $S\subseteq V(G) \ (|S|\geq 2)$\\[0.1cm]
\cline{1-3} Set of Steiner trees & $\left\{
\begin{array}{ll}
\mathscr{P}_{x,y}=\{P_1,P_2,\cdots,P_{\ell}\}\\[0.1cm]
\{x,y\}\subseteq V(P_i),\\[0.1cm]
E(P_i)\cap
E(P_j)=\varnothing\\[0.1cm]
V(P_i)\cap V(P_j)=\{x,y\}\\[0.1cm]
\end{array}
\right.$ & $\left\{
\begin{array}{ll}
\mathscr{T}_{S}=\{T_1,T_2,\cdots,T_{\ell}\}\\[0.1cm]
S\subseteq V(T_i),\\[0.1cm]
E(T_i)\cap E(T_j)=\varnothing,\\[0.1cm]
V(T_i)\cap
V(T_j)=S\\[0.1cm]
\end{array}
\right.$\\[0.1cm]
\cline{1-3}
Local parameter & $\kappa(x,y)=\max|\mathscr{P}_{x,y}|$ & $\kappa(S)=\max|\mathscr{T}_{S}|$\\[0.1cm]
\cline{1-3} Global parameter & $\kappa(G)=\underset{x,y\in
V(G)}{\min} \kappa(x,y)$ & $\kappa_k(G)=\underset{S\subseteq V(G),
|S|=k}{\min} \kappa(S)$\\[0.1cm]
\cline{1-3}
\end{tabular}
\end{center}
\begin{center}
{Table 1. Classical~connectivity and generalized connectivity}
\end{center}
}

As a natural counterpart of the generalized connectivity, we
introduced in \cite{LMS} the concept of generalized
edge-connectivity, which is a generalization of the `path' version
definition of edge-connectivity. For $S\subseteq V(G)$ and $|S|\geq
2$, the \emph{generalized local edge-connectivity} $\lambda(S)$ is
the maximum number of edge-disjoint Steiner trees connecting $S$ in
$G$. For an integer $k$ with $2\leq k\leq n$, the \emph{generalized
$k$-edge-connectivity} $\lambda_k(G)$ of $G$ is then defined as
$\lambda_k(G)= \min\{\lambda(S)\,|\,S\subseteq V(G) \ and \
|S|=k\}$. It is also clear that when $|S|=2$, $\lambda_2(G)$ is
nothing new but the standard edge-connectivity $\lambda(G)$ of $G$,
that is, $\lambda_2(G)=\lambda(G)$, which is the reason why we
address $\lambda_k(G)$ as the generalized edge-connectivity of $G$.
Also set $\lambda_k(G)=0$ when $G$ is disconnected. Results on the
generalized edge-connectivity can be found in \cite{LM2, LMS, LMS}.

\subsection{Path connectivity and path edge-connectivity}

Dirac \cite{Dirac} showed that in a $(k-1)$-connected graph there is
a path through each $k$ vertices. Related problems were inquired in
\cite{Wilson}. In \cite{Hager2}, Hager revised this statement to the
question of how many internally disjoint paths $P_i$ with the
exception of a given set $S$ of $k$ vertices exist such that
$S\subseteq V(P_i)$. The path connectivity of a graph $G$,
introduced by Hager \cite{Hager2}, is a natural specialization of
the generalized connectivity and is also a natural generalization of
the `path' version definition of connectivity. For a graph $G=(V,E)$
and a set $S\subseteq V(G)$ of at least two vertices, \emph{a path
connecting $S$} (or simply, \emph{an $S$-path}) is a subgraph
$P=(V',E')$ of $G$ that is a path with $S\subseteq V'$. Note that a
path connecting $S$ is also a tree connecting $S$. Two paths $P$ and
$P'$ connecting $S$ are said to be \emph{internally disjoint} if
$E(P)\cap E(P')=\varnothing$ and $V(P)\cap V(P')=S$. For $S\subseteq
V(G)$ and $|S|\geq 2$, the \emph{local path connectivity} $\pi_G(S)$
is the maximum number of internally disjoint paths connecting $S$ in
$G$, that is, we search for the maximum cardinality of edge-disjoint
paths which contain $S$ and are vertex-disjoint with the exception
of the vertices in $S$. For an integer $k$ with $2\leq k\leq n$, the
\emph{path $k$-connectivity} is defined as
$\pi_k(G)=\min\{\pi_G(S)\,|\,S\subseteq V(G),|S|=k\}$, that is,
$\pi_k(G)$ is the minimum value of $\pi_G(S)$ when $S$ runs over all
$k$-subsets of $V(G)$. Clearly, $\pi_1(G)=\delta(G)$ and
$\pi_2(G)=\kappa(G)$. For $k\geq 3$, $\pi_k(G)\leq \kappa_k(G)$
holds because each path is also a tree.

The relation between generalized connectivity and path connectivity
are shown in the following Table 2. {\small
\begin{center}
\begin{tabular}{|c|c|c|}
\hline  & Generalized connectivity&
Path-connectivity\\[0.1cm]
\cline{1-3}
Vertex subset & $S\subseteq V(G) \ (|S|\geq 2)$ & $S\subseteq V(G) \ (|S|\geq 2)$\\[0.1cm]
\cline{1-3} Set of Steiner trees & $\left\{
\begin{array}{ll}
\mathscr{T}_{S}=\{T_1,T_2,\cdots,T_{\ell}\}\\
S\subseteq V(T_i),\\
E(T_i)\cap E(T_j)=\varnothing,\\
\end{array}
\right.$ & $\left\{
\begin{array}{ll}
\mathscr{P}_{S}=\{P_1,P_2,\cdots,P_{\ell}\}\\
S\subseteq V(P_i),\\
E(P_i)\cap E(P_j)=\varnothing,\\
\end{array}
\right.$\\[0.05cm]
\cline{1-3}
Local parameter  & $\kappa(S)=\max|\mathscr{T}_{S}|$ & $\pi(S)=\max|\mathscr{P}_{S}|$\\
\cline{1-3} Global parameter & $\kappa_k(G)=\underset{S\subseteq
V(G), |S|=k}{\min} \kappa(S)$ & $\pi_k(G)=\underset{S\subseteq V(G),
|S|=k}{\min} \pi(S)$\\[0.05cm]
\cline{1-3}
\end{tabular}
\end{center}
\begin{center}
{Table 2. Two kinds of tree-connectivities}
\end{center}
}

As a natural counterpart of path $k$-connectivity, we recently
introduced the concept of path $k$-edge-connectivity. Two paths $P$
and $P'$ connecting $S$ are said to be \emph{edge-disjoint} if
$E(P)\cap E(P')=\varnothing$. For $S\subseteq V(G)$ and $|S|\geq 2$,
the \emph{path local edge-connectivity} $\omega_G(S)$ is the maximum
number of edge-disjoint paths connecting $S$ in $G$. For an integer
$k$ with $2\leq k\leq n$, the \emph{path $k$-edge-connectivity} is
defined as $\omega_k(G)=\min\{\omega_G(S)\,|\,S\subseteq
V(G),|S|=k\}$, that is, $\omega_k(G)$ is the minimum value of
$\omega_G(S)$ when $S$ runs over all $k$-subsets of $V(G)$. Clearly,
we have
$$
\ \ \ \ \ \ \ \ \ \ \ \ \ \ \ \ \left\{
\begin{array}{ll}
\omega_k(G)=\delta(G),&\mbox {\rm for}~k=1;\\
\omega_k(G)=\lambda(G),&\mbox {\rm for}~k=2;\ \ \ \ \ \ \ \ \ \ \  \
\ \ \ \ \ \ \ \ \ \ \ \ \ \ \ \ \ \ \ \ \ \ \ \ \ \ \ \
\ (1)\\
\omega_k(G)\leq \lambda_k(G),&\mbox {\rm for} \ k\geq 3.
\end{array}
\right.
$$

The following observation is immediate.
\begin{obs}\label{obs1-1}
$(1)$ Let $G$ be a connected graph. Then $\pi_k(G)\leq
\omega_k(G)\leq \delta(G)$.

$(2)$ Let $G$ be a connected graph with minimum degree $\delta$. If
$G$ has two adjacent vertices of degree $\delta$, then $\pi_k(G)\leq
\omega_k(G)\leq \delta-1$.
\end{obs}

\begin{lem}{\upshape\cite{Hager}}\label{lem1-2}
Let $k,n$ be two integers with $3\leq k\leq n$, and let $K_n$ be a
complete graph of order $n$. Then
$$
\pi_k(K_n)=\left\lfloor\frac{2n+k^2-3k}{2(k-1)}\right\rfloor.
$$
\end{lem}

Note that each graph is a spanning subgraph of a complete graph. So
the following result is immediate.
\begin{obs}\label{obs1-3}
Let $k,n$ be two integers with $3\leq k\leq n$, and let $G$ be a
graph of order $n$. Then
$$
0\leq \pi_k(G)\leq
\left\lfloor\frac{2n+k^2-3k}{2(k-1)}\right\rfloor.
$$
\end{obs}

\begin{lem}\label{lem1-4}
Let $G$ be a graph of order $n$. If $n$ is even, then
$\pi_3(G)=\frac{n}{2}$ if and only if $G$ is a complete graph of
order $n$.
\end{lem}
\begin{pf}
From Lemma \ref{lem1-2}, we have $\pi_3(K_n)=\lfloor
\frac{n}{2}\rfloor=\frac{n}{2}$. Actually, the complete graph $K_n$
is the unique graph with this property. We only need to show that
$\pi_k(K_n\setminus uv)<\frac{n}{2}-1$ for any $uv\in E(K_n)$. Pick
one vertex $w\in V(K_n)-\{u,v\}$. Choose $S=\{u,v,w\}$. If we choose
$uvw$ as a path, then any other $S$-Steiner path occupies at least
two vertices of $V(G)-S$, and hence $\pi(S)\leq
1+\lfloor\frac{n-3}{2}\rfloor=\frac{n}{2}-1$. Suppose that $uv$ and
$vw$ belong to different $S$-Steiner path. Then $\pi(S)\leq
2+\lfloor\frac{n-5}{2}\rfloor=\frac{n}{2}-1$. So $\pi_3(G)\leq
\pi(S)\leq \frac{n}{2}-1$, as desired.\qed
\end{pf}
\begin{lem}{\upshape \cite{Hager2}}\label{lem1-5}
For any connected graph $G$, $\pi_3(G)\geq \frac{1}{2}\kappa(G)$.
Moreover, the lower bound is sharp.
\end{lem}

In \cite{Jain}, Jain et al. obtained the following result.
\begin{lem}{\upshape\cite{Jain}}\label{lem1-6}
Let $G$ be a graph and $S=\{v_1,v_2,v_3\}$ be a subset of $V(G)$.
Assume that $v_1$ and $v_2$ are $\ell$-edge-connected, and $v_1$ and
$v_3$ are $\frac{\ell}{2}$-edge-connected in $G$. Then $G$ has
$\frac{\ell}{2}$ edge-disjoint $S$-Steiner trees.
\end{lem}

In their proof, the $\frac{\ell}{2}$ edge-disjoint $S$-Steiner trees
are actually edge-disjoint $S$-Steiner paths. So the following
result is immediate.
\begin{cor}\label{cor1-5}
For any connected graph $G$, $\omega_3(G)\geq
\frac{1}{2}\lambda(G)$. Moreover, the lower bound is sharp.
\end{cor}

\subsection{Application background of these parameters}

In addition to being a natural combinatorial measure, the path
$k$-(edge-)connectivity and generalized $k$-(edge-)connectivity can
be motivated by their interesting interpretation in practice. For
example, suppose that $G$ represents a network. If one considers to
connect a pair of vertices of $G$, then a path is used to connect
them. However, if one wants to connect a set $S$ of vertices of $G$
with $|S|\geq 3$, then a tree has to be used to connect them. This
kind of tree for connecting a set of vertices is usually called a
{\it Steiner tree}, and popularly used in the physical design of
VLSI circuits (see \cite{Grotschel1, Grotschel2, Sherwani}). In this
application, a Steiner tree is needed to share an electric signal by
a set of terminal nodes. Steiner tree is also used in computer
communication networks (see \cite{Du}) and optical wireless
communication networks (see \cite{Cheng}). Usually, one wants to
consider how tough a network can be, for the connection of a set of
vertices. Then, the number of totally independent ways to connect
them is a measure for this purpose. The path $k$-connectivity and
generalized $k$-connectivity can serve for measuring the capability
of a network $G$ to connect any $k$ vertices in $G$.

\subsection{Main results of this paper}

Chartrand and Stewart \cite{Steeart} investigated the relation
between the connectivity and edge-connectivity of a graph and its
line graph.
\begin{thm}{\upshape \cite{Steeart}}\label{th1-8}
If $G$ is a connected graph, then

$(1)$ $\kappa(L(G))\geq \lambda(G)$ if $\lambda(G)\geq 2$.

$(2)$ $\lambda(L(G))\geq 2\lambda(G)-2$.

$(3)$ $\kappa(L(L(G)))\geq 2\kappa(G)-2$.
\end{thm}

In Section 2, we investigate the relation between the path
connectivity and path edge-connectivity of a graph and its line
graph.

In their book, Capobianco and Molluzzo \cite{Capobianco}, using
$K_{1,n}$ as their example, note that the difference between the
connectivity of a graph and its line graph can be arbitrarily large.
They then proposed an open problem: Whether for any two integers
$p,q \ (1<p<q)$, there exists a graph $G$ such that $\kappa(G)=p$
and $\kappa(L(G))=q$. In \cite{Bauer}, Bauer and Tindell gave a
positive answer of this problem, that is, for every pair of integers
$p,q \ (1<p<q)$ there is a graph of connectivity $p$ whose line
graph has connectivity $q$.

Note that the difference between the path $k$-connectivity of a
graph $G$ and its line graph $L(G)$ can be also arbitrarily large.
Let $n,k$ be two integers with $2\leq k\leq n$, and let $G=K_{1,n}$.
Then $L(G)=K_n$, $\pi_k(G)=0$ and
$\pi_k(L(G))=\left\lfloor\frac{2n+k^2-3k}{2(k-1)}\right\rfloor$. In
fact, we can consider a similar problem: Whether for any two
integers $p,q$, $1<p<q$, there exists a graph $G$ such that
$\pi_k(G)=p$ and $\pi_k(L(G))=q$. It seem to be not easy to solve
this problem for a general $k$. In this paper, we focus our
attention on the case $k=3$ and $q\geq 2p-1$, and give a positive
answer of this problem.

\section{Path connectivity of line graphs}

Hager \cite{Hager} obtain the following result.
\begin{lem}{\upshape\cite{Hager}}\label{lem2-1}
Let $G$ be a graph of order $n$, and let $k$ be an integer with
$3\leq k\leq n$. If $\kappa(G)\geq 2^{k-2}\ell$, then $\pi_k(G)\geq
\ell$.
\end{lem}
\begin{cor}\label{cor2-2}
Let $G$ be a graph of order $n$, and let $k$ be an integer with
$3\leq k\leq n$. Then
$$
2^{2-k}\kappa(G)\leq \pi_k(G)\leq \kappa(G).
$$
\end{cor}

For general $k$, we can easily derive the following result.
\begin{pro}\label{th2-3}
Let $G$ be a graph of order $n$. If $G$ is $2$-edge-connected, then

$(1)$ $2^{2-k}\omega_{k}(G)\leq \pi_{k}(L(G))$;

$(2)$ $\omega_{k}(L(G))\geq 2^{2-k}\omega_k(G)$;

$(3)$ $\pi_{k}(L(L(G))\geq 2^{3-k}(\pi_k(G)-1)$.
\end{pro}
\begin{pf}
For $(1)$, from Corollary \ref{cor2-2} and Theorem \ref{th1-8}, we
have
$$
\pi_k(L(G))\geq 2^{2-k}\kappa(L(G)) \geq 2^{2-k}\lambda(G)\geq
2^{2-k}\kappa(G)\geq 2^{2-k}\pi_k(G).
$$

For $(2)$, from $(1)$ of this proposition, Corollary \ref{cor2-2}
and Theorem \ref{th1-8}, we have
$$
\omega_k(L(G))\geq \pi_k(L(G))\geq 2^{2-k}\kappa(L(G)) \geq
2^{2-k}\lambda(G)\geq 2^{2-k}\omega_k(G).
$$

For $(3)$, from $(1)$ of this proposition and Corollary
\ref{cor2-2}, we have
$$
\pi_k(L(L(G)))\geq 2^{2-k}\kappa(L(L(G)))\geq 2^{2-k}(2\kappa(G)-2)=2^{3-k}(\kappa(G)-1)\\
\geq 2^{3-k}(\pi_k(G)-1).
$$
\end{pf}

For $k=3$, we can improve the above result, and obtain the following
theorem.
\begin{thm}\label{th2-4}
Let $G$ be a $2$-connected. Then

$(1)$ $\omega_3(G)\leq \pi_3(L(G))$.

$(2)$ $\omega_3(L(G))\geq \omega_3(G)-1$.

$(3)$ $\pi_3(L(L(G))\geq \pi_3(G)-1$.
\end{thm}
\begin{pf}
For $(1)$, let $e_1,e_2,e_3$ be three arbitrary distinct vertices of
the line graph of $G$ such that $\omega_3(G)=\ell$ with $\ell \geq
1$. Let $e_1=v_1v_1'$, $e_2=v_2v_2'$ and $e_3=v_3v_3'$ be those
edges of $G$ corresponding to the vertices $e_1,e_2,e_3$ in $L(G)$,
respectively.

Consider three distinct vertices of the six end-vertices of
$e_1,e_2,e_3$. Without loss of generality, let $S=\{v_1,v_2,v_3\}$
be three distinct vertices. Since $\omega_3(G)=\ell$, there exist
$\ell$ edge-disjoint $S$-Steiner paths in $G$, say
$P_1,P_2,\cdots,P_{\ell}$. We define a minimal $S$-Steiner path $P$
as an $S$-Steiner path whose sub-path obtained by deleting any edge
of $P$ does not connect $S$.
\begin{figure}[h,t,b,p]
\begin{center}
\scalebox{0.7}[0.7]{\includegraphics{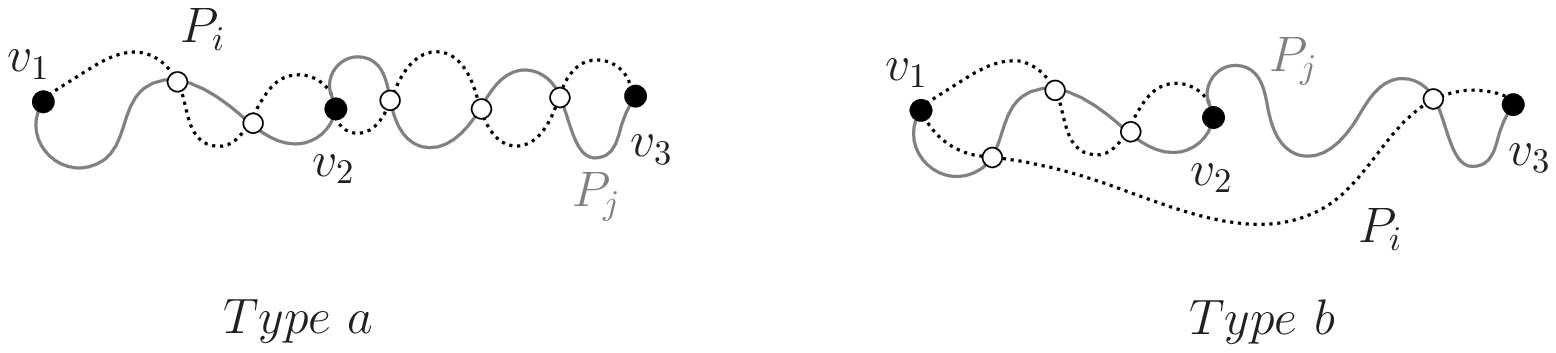}}\\
Figure 1:  Three possible types of $T_i\cup T_j$.
\end{center}
\end{figure}
Choosing any two edge-disjoint minimal $S$-Steiner paths $P_i$ and
$P_j \ (1\leq i,j\leq \ell)$ in $G$, we will show that the paths
$P_i'$ and $P_j'$ corresponding to $P_i$ and $P_j$ in $L(G)$ are
internally disjoint $S$-Steiner paths. It is easy to see that
$P_i\cup P_j$ has two possible types, as shown in Figure $1$. Since
$P_i$ and $P_j$ are edge-disjoint in $G$, we can find internally
disjoint Steiner paths $P_i'$ and $P_j'$ connecting $e_1,e_2,e_3$ in
$L(G)$. We give an example of Type $a$, see Figure $2$. So
$\pi_3(L(G))\geq \ell$, as desired.
\begin{figure}[h,t,b,p]
\begin{center}
\scalebox{0.6}[0.6]{\includegraphics{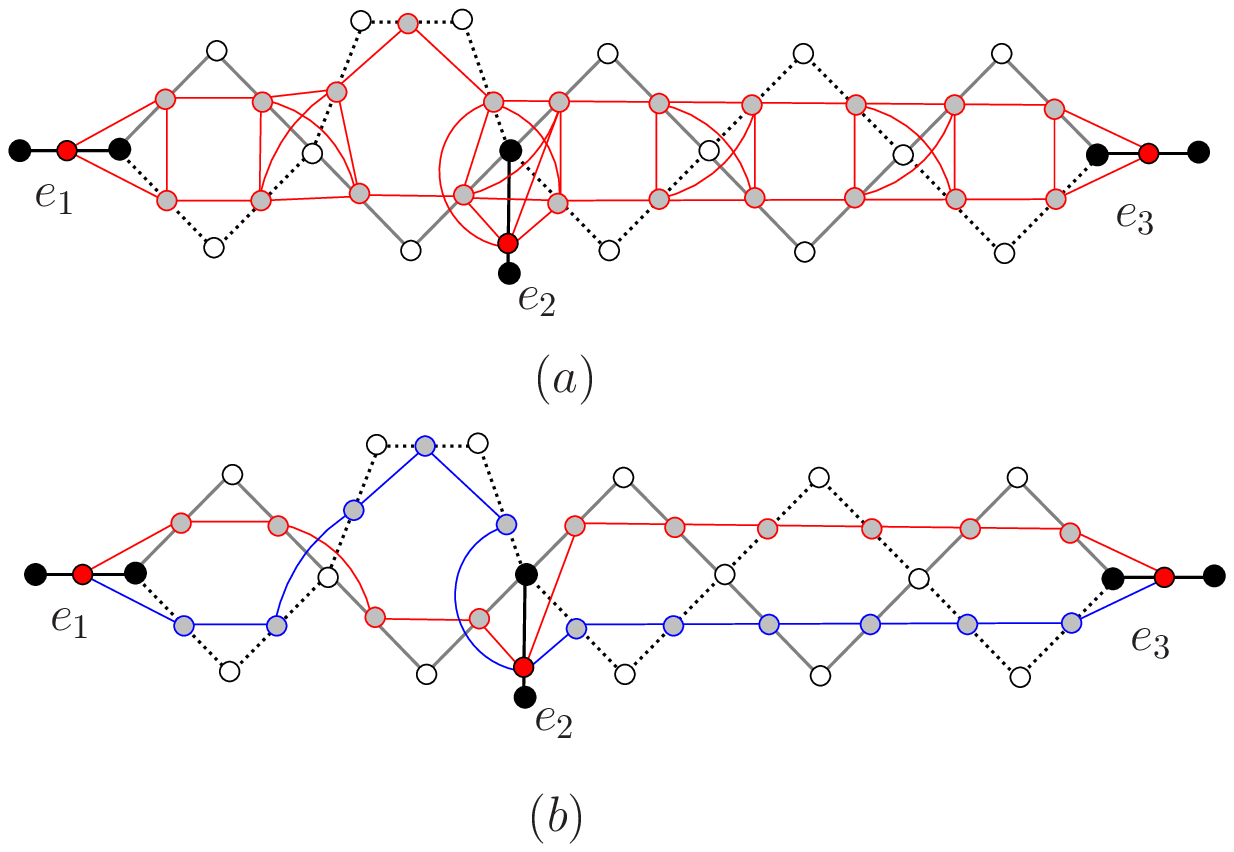}}\\
Figure 2 (a): An example for $P_i$ and $P_j$ connecting $S$ and their line graphs. \\
Figure 2 (b): An example for $P_i'$ and $P_j'$ corresponding to
$P_i$ and $P_j$.
\end{center}
\end{figure}

$(2)$ From Corollary \ref{cor1-5} and Theorem \ref{th1-8}, we have
$$
\omega_3(L(G))\geq \frac{1}{2}\lambda(L(G)) \geq
\frac{1}{2}(2\lambda(G)-2)=\lambda(G)-1\geq \omega_3(G)-1.
$$

$(3)$ From Lemma \ref{lem1-5} and Theorem \ref{th1-8}, we have
$$
\pi_3(L(L(G)))\geq \frac{1}{2}\kappa(L(L(G)))
\geq\frac{1}{2}(2\kappa(G)-2)=\kappa(G)-1\geq \pi_3(G)-1.
$$
\end{pf}

\section{Graphs with prescribed
path connectivity and path edge-connectivity}

In {\upshape\cite{Hager}}, Hager got the following result.
\begin{lem}{\upshape\cite{Hager}}\label{lem3-1}
Let $K_{a,b}$ be a complete bipartite graph with $a+b$ vertices, and
let $k$ be an integer with $2\leq k\leq a+b$. Then
$$
\pi_k(K_{a,b})=\min\left\{\frac{a}{k-1}, \frac{b}{k-1}\right\}.
$$
\end{lem}

The following corollary is immediate from the above theorem.
\begin{cor}\label{cor3-2}
Let $a,b$ be two integers with $2\leq a \leq b$, and $K_{a,b}$
denote a complete bipartite graph with a bipartition of sizes $a$
and $b$, respectively. Then
$$
\pi_3(K_{a,b})=\left \lfloor \frac{a}{2}\right \rfloor.
$$
\end{cor}

In this section, we consider the problem mentioned in Subsection 1.5
for the case $q\geq 2p-1$. Let us put our attention on the complete
bipartite graph $G=K_{2p,2q-2p+2}$. Since $q\geq 2p-1$, it follows
that $2q-2p+2\geq 2p$. From Corollary \ref{cor3-2},
$\pi_3(G)=\pi_3(K_{2p,2q-2p+2})=p$. Now we turn to consider the line
graph of complete bipartite graph $G=K_{2p,2q-2p+2}$.

Recall that the Cartesian product (also called the {\em square
product}) of two graphs $G$ and $H$, written as $G\Box H$, is the
graph with vertex set $V(G)\times V(H)$, in which two vertices
$(u,v)$ and $(u',v')$ are adjacent if and only if $u=u'$ and
$(v,v')\in E(H)$, or $v=v'$ and $(u,u')\in E(G)$. Clearly, the
Cartesian product is commutative, that is, $G\Box H\cong H\Box G$.

The following lemma is from \cite{West}.
\begin{lem}{\upshape \cite{West}}\label{lem3-3}
For a complete bipartite graph $K_{r,s}$, $L(K_{r,s})=K_r\Box K_s$.
\end{lem}

From the above lemma, $L(G)=L(K_{2p,2q-2p+2})=K_{2p}\Box
K_{2q-2p+2}$. In order to obtain the exact value of
$\pi_3(L(G))=\pi_3(K_{2p}\Box K_{2q-2p+2})$, we consider to
determine the exact value of the Cartesian product of two complete
graphs.

Before proving Lemma \ref{lem3-4}, we define some notation. Let $G$
and $H$ be two connected graphs with $V(G)=\{u_1,u_2,\ldots,u_{n}\}$
and $V(H)=\{v_1,v_2,\ldots,v_{m}\}$, respectively. Then $V(G\circ
H)=\{(u_i,v_j)\,|\,1\leq i\leq n, \ 1\leq j\leq m\}$. For $v\in
V(H)$, we use $G(v)$ to denote the subgraph of $G\circ H$ induced by
the vertex set $\{(u_i,v)\,|\,1\leq i\leq n\}$. Similarly, for $u\in
V(G)$, we use $H(u)$ to denote the subgraph of $G\circ H$ induced by
the vertex set $\{(u,v_j)\,|\,1\leq j\leq m\}$.
\begin{lem}\label{lem3-4}
Let $p,q$ be two integers with $p\geq 2$ and $q\geq 3$. Then
$$
\pi_3(K_{2p}\Box K_{2q-2p+2})=q.
$$
\end{lem}
\begin{pf}
Let $G=K_{2p}$ and $H=K_{2q-2p+2}$. Set
$V(G)=\{u_1,u_2,\cdots,u_{2p}\}$ and
$V(H)=\{v_1,v_2,\cdots,v_{2q-2p+2}\}$. From Lemma \ref{lem1-4}, we
have
$$
\pi_3(G\Box H)=\pi_3(K_{2p}\Box K_{2q-2p+2})\leq \pi_3(K_{2q+2})\leq
q.
$$
It suffices to show that $\pi_3(G\Box H)\geq q$. We need to show
that for any $S=\{x,y,z\}\subseteq V(G\Box H)$, there exist $q$
internally disjoint $S$-Steiner paths. We complete our proof by the
following three cases.

\textbf{Case 1}. $x,y,z$ belongs to the same $V(H(u_i)) \ (1\leq
i\leq r)$.

Without loss of generality, we assume $x,y,z\in V(H(u_1))$. Since
$\pi_3(H)=\pi_3(K_{2q-2p+2})$ $=q-p+1$, there exist $q-p+1$
internally disjoint $S$-Steiner paths $P_1,P_2,\cdots,P_{q-p+1}$ in
$H(u_1)$. Let $x_j,y_j,z_j$ be the vertices corresponding to $x,y,z$
in $H(u_j) \ (2\leq j\leq 2p)$. Then the paths $Q_{i}$ induced by
the edges in
$\{xx_{2i},x_{2i}y_{2i},y_{2i}y,yy_{2i+1},y_{2i+1}z_{2i+1},z_{2i+1}z\}
\ (1\leq i\leq p-1)$ are $p-1$ internally disjoint $S$-Steiner
paths. These paths together with the paths
$P_1,P_2,\cdots,P_{q-p+1}$ are $q$ internally disjoint $S$-Steiner
paths, as deisred.

\textbf{Case 2}. Only two vertices of $\{x,y,z\}$ belong to some
copy $H(u_i)$.

We may assume $x,y\in V(H(u_1))$, $z\in V(H(u_2))$. Let $x',y'$ be
the vertices corresponding to $x,y$ in $H(u_2)$, and let $z'$ be the
vertex corresponding to $z$ in $H(u_1)$.

Suppose $z'\not\in \{x,y\}$. Without loss of generality, let
$x=(u_1,v_{1})$, $y=(u_1,v_2)$ and $z'=(u_1,v_3)$. Then the path
$Q_1$ induced by the edges in $\{xy,yy',y'z\}$, the path $Q_2$
induced by the edges in $\{xx',x'z,zz',z'y\}$ and $P_{i}$ induced by
the edges in $\{x(u_1,v_{2i}),(u_1,v_{2i})y,$
$y(u_1,v_{2i+1}),(u_1,v_{2i+1})(u_2,v_{2i+1}),(u_2,v_{2i+1})z\} \
(2\leq i\leq q-p)$ are $q-p+1$ internally disjoint $S$-Steiner
paths. Note that all the edges of these paths are between $H(u_1)$
and $H(u_2)$. Let $x_j,y_j,z_j$ be the vertices corresponding to
$x,y,z'$ in $H(u_j) \ (3\leq j\leq 2p)$. The the paths $Q_{j}$
induced by the edges in
$\{xx_{2j+1},x_{2j+1}y_{2j+1},x_{2j+1}y,yy_{2j+2},y_{2j+2}z_{2j+2},z_{2j+2}z\}
\ (1\leq j\leq p-1)$ are $p-1$ internally disjoint $S$-Steiner
paths. These paths together with the paths
$Q_1,Q_2,P_2,P_3,\cdots,P_{q-p}$ are $q$ internally disjoint
$S$-Steiner paths, as desired.

Suppose $z'\in \{x,y\}$. Without loss of generality, let
$x=(u_1,v_{1})$ and $y=(u_1,v_2)$. Then the path $Q$ induced by the
edges in $\{xy,yz\}$ and $P_{j}$ induced by the edges in
$\{x(u_1,v_{2j-1}),(u_1,v_{2j-1})y,y(u_1,v_{2j}),(u_1,v_{2j})(u_2,v_{2j}),(u_2,v_{2j})z\}
\ (2\leq j\leq q-p+1)$ are $q-p+1$ internally disjoint $S$-Steiner
paths. Let $x_j,y_j,z_j$ be the vertices corresponding to $x,y,z'$
in $H(u_j) \ (3\leq j\leq 2p)$. Then the paths $Q_{i}$ induced by
the edges in
$\{xx_{2i+1},x_{2i+1}y_{2i+1},y_{2i+1}y,yy_{2i+2},y_{2i+2}z_{2i+2},z_{2i+2}z\}
\ (1\leq i\leq p-1)$ are $p-1$ internally disjoint $S$-Steiner
paths. These paths together with the paths
$Q,P_1,P_2,\cdots,P_{q-p}$ are $q$ internally disjoint $S$-Steiner
paths, as desired.

\textbf{Case 3}. $x,y,z$ are contained in distinct $H(u_i)$s.

We may assume that $x\in V(H(u_1))$, $y\in V(H(u_2))$, $z\in
V(H(u_3))$. Let $y',z'$ be the vertices corresponding to $y,z$ in
$H(u_1)$, $x',z''$ be the vertices corresponding to $x,z$ in
$H(u_2)$ and $x'',y''$ be the vertices corresponding to $x,y$ in
$H(u_3)$.

Suppose that $x,y',z'$ are distinct vertices in $H(u_1)$. Without
loss of generality, let $x=(u_1,v_1)$, $y=(u_2,v_2)$ and
$z=(u_3,v_3)$. Let $x_j,y_j,z_j$ be the vertices corresponding to
$x,y,z'$ in $H(u_j) \ (4\leq j\leq 2p)$. Then the path $R_1$ induced
by the edges in $\{xx',yx',yy'',y''z\}$, the path $R_2$ induced by
the edges in $\{xy',yy',yz'',z''z\}$, the path $R_3$ induced by the
edges in $\{xz',zz',yy_{2p},y_{2p}z_{2p},z_{2p}z\}$ and $P_{j}$
induced by the edges in
$\{x(u_1,v_{2j}),(u_1,v_{2j})(u_2,v_{2j}),y(u_2,v_{2j}),y(u_2,v_{2j+1}),(u_2,v_{2j+1})(u_3,v_{2j+1}),(u_3,v_{2j+1})z\}
\ (2\leq j\leq q-p)$ are $q-p+2$ internally disjoint $S$-Steiner
paths. The the paths $Q_{i}$ induced by the edges in
$\{xx_{2i-1},x_{2i-1}y_{2i-1},y_{2i-1}y,yy_{2i},y_{2i}z_{2i},z_{2i}z\}
\ (2\leq i\leq p-1)$ are $p-2$ internally disjoint $S$-Steiner
paths. These paths together with the trees
$R_1,R_2,R_3,P_2,P_3,\cdots,P_{q-p-1}$ are $q$ internally disjoint
$S$-Steiner paths, as desired.

Suppose that two of $x, y',z'$ are the same vertex in $H(u_1)$.
Without loss of generality, let $x=(u_1,v_1)$, $y=(u_2,v_1)$ and
$z=(u_3,v_2)$. Let $x_j,y_j,z_j$ be the vertices corresponding to
$x,y,z'$ in $H(u_j) \ (4\leq j\leq 2p)$. Then the path $R_1$ induced
by the edges in $\{xy,yx'',x''z\}$, the path $R_2$ induced by the
edges in $\{xx_{2p},x_{2p}y_{2p},y_{2p}y,yz'',z''z\}$ and $P_{j}$
induced by the edges in
$\{x(u_1,v_{2j-1}),(u_1,v_{2j-1})(u_2,v_{2j-1}),(u_2,v_{2j-1})y,y(u_2,v_{2j}),(u_2,v_{2j})(u_3,v_{2j}),(u_3,v_{2j})z\}
$\\ $(2\leq j\leq q-p+1)$ are $q-p+2$ internally disjoint
$S$-Steiner paths. The the paths $Q_{i}$ induced by the edges in
$\{xx_{2i-1},x_{2i-1}y_{2i-1},y_{2i-1}y,yy_{2i},y_{2i}z_{2i},z_{2i}z\}
\ (2\leq i\leq p-1)$ are $p-2$ internally disjoint $S$-Steiner
paths. These paths together with the paths
$R_1,R_2,P_2,P_3,\cdots,P_{q-p+1}$ are $q$ internally disjoint
$S$-Steiner paths, as desired.

Suppose that $x,y',z'$ are the same vertex in $H(u_1)$. Without loss
of generality, let $x=(u_1,v_1)$, $y=(u_2,v_1)$ and $z=(u_3,v_1)$.
Then the path $R_1$ induced by the edges in $\{xy,yz\}$, the path
$R_2$ induced by the edges in
$\{xz,y(u_2,v_{2q-2p+2}),(u_2,v_{2q-2p+2})$
$(u_3,v_{2q-2p+2}),(u_3,v_{2q-2p+2})z\}$ and $P_{j}$ induced by the
edges in $\{x(u_1,v_{2j}),(u_1,v_{2j})(u_2,v_{2j}),$
$y(u_2,v_{2j}),y(u_2,v_{2j+1}),(u_2,v_{2j+1})(u_3,v_{2j+1}),(u_3,v_{2j+1})z\}
\ (1\leq j\leq q-p)$ are $q-p+2$ internally disjoint $S$-Steiner
paths. Let $x_j,y_j,z_j$ be the vertices corresponding to $x,y,z'$
in $H(u_j) \ (4\leq j\leq 2p)$. The the paths $Q_{i}$ induced by the
edges in $\{xx_{2i},x_{2i}y_{2i},y_{2i}y,yy_{2i+1},$
$y_{2i+1}z_{2i+1} ,z_{2i+1}z\} \ (2\leq i\leq p-1)$ are $p-2$
internally disjoint $S$-Steiner paths. These paths together with the
paths $R_1,R_2,P_1,P_2,\cdots,P_{q-p}$ are $q$ internally disjoint
$S$-Steiner paths, as desired.

From the above argument, we conclude that for any
$S=\{x,y,z\}\subseteq V(G\Box H)$, there exist $q$ internally
disjoint $S$-Steiner paths, and hence $\pi(S)\geq q$. From the
arbitraries of $S$, we have $\pi_k(G)=q$. \qed
\end{pf}

From Lemma \ref{lem3-4}, the following result holds.
\begin{thm}\label{pro3-3}
For any two integers $p,q$ with $q\geq 2p-1$, $p\geq 2$ and $q\geq
3$, there exists a graph $G$ such that $\pi_3(G)=p$ and
$\pi_3(L(G))=q$.
\end{thm}

\end{document}